\documentclass{amsart} 
\usepackage{amsfonts,amsmath,amsthm,amssymb}
\begin{document}

\title{On the isomorphism problem for unit groups of modular group algebras}
\author{A.~Konovalov, A.~ Krivokhata}
\date{}

\subjclass{MSC 16S34, 20C05}

\thanks{The first author was supported by Francqui Stichting (Belgium) grant ADSI107}

\begin{abstract}
Using the computational algebra system \textsf{GAP} and the 
package \textsf{LAGUNA}, we checked that all
2-groups of order not greater than 32 are determined 
by normalized unit groups of their modular group algebras
over the field of two elements.
\end{abstract}

\maketitle

\section{Introduction}
\label{Intro}

Let $KG$ be a modular group algebra of a finite $p$-group $G$ over
the field $K$ of $p$ elements, and $V(KG)$ be the normalized unit
group of $KG$. 

The Modular Isomorphism Problem ({\bf MIP}) asks whether a finite
$p$-group $G$ is determined by its modular group algebra over the field
of $p$ elements, i.e. is it true that
$$
KG \cong KH \Rightarrow G \cong H \; ?
$$
It remains open for more than 50 years, and up to now it is solved
only for some classes of groups and for some groups of small orders.

In the context of {\bf MIP}, S.~D.~Berman formulated 
the Modular Isomorphism Problem for Normalized Unit Groups ({\bf UMIP}),
asking whether a finite $p$-group $G$ is determined by the normalized 
unit group of its modular group algebra over the field of $p$ elements, i.e.
is it true that
$$
V(KG) \cong V(KH) \Rightarrow G \cong H \; ?
$$

Obviously, {\bf UMIP} is stronger than {\bf MIP}, and from the solution
of {\bf UMIP} follows the solution of {\bf MIP}.
For a long time, the positive solution of {\bf UMIP} 
was known only for abelian $p$-groups, and only recently it was solved for 
2-groups of maximal class in \cite{BB1}, and for $p$-groups with the cyclic
Frattini subgroup for $p>2$ in \cite{BB2}.

Although {\bf UMIP} is a stronger problem, it is well-suited for computer-aided
tests. Since the order of the group $G$ is determined by $V(KG)$, 
for such a test it will be enough to check that all non-isomorphic 
groups of a given order have non-isomorphic normalized unit groups 
of their modular group algebras over the
field of 2 elements. The structure of the normalized unit group $V(KG)$ can be
investigated using the computational algebra system \textsf{GAP} 
\cite{GAP} enhanced by the \textsf{LAGUNA} package \cite{LAGUNA} to calculate
the power-commutator presentation of $V(KG)$. Having such a presentation, we can compare various invariants of normalized unit groups such as their exponent, sizes of certain subgroups, the number of involutions, etc. until we find a pair of different ones making evidence of their non-isomorphism. Note that here we have much more freedom in choosing invariants unlike in the similar approach to {\bf MIP} tests when we can compare only properties of $G$ known to be determined by $KG$.

In the present paper we give a detailed report about checking that all 2-groups of order not greater than 32 are determined by the normalized unit group of their modular group algebras. Thus, our tests confirm {\bf UMIP} for these groups.
 
\section{Groups of orders up to 16}
\label{order16}

Since the solution of {\bf UMIP} is known for abelian groups and also for
2-groups of maximal class \cite{BB1}, there is no need to consider
groups of orders 2, 4, and 8, so we can start from order 16. 

Among 14 groups of order 16, there are 
5 abelian groups and 3 groups of maximal class, so it remains to consider
6 groups of this order.

In the following table for these groups we give the number of the group $G$ 
in the library of groups of order 16 from the Small Groups Library of the 
\textsf{GAP} system, the order of the Frattini subgroup of $V(KG)$, and the 
number of elements of order 2 in $V(KG)$:

$$
\begin{array}{|c|c|c|}\hline
\text{Catalogue} & \text{Order of the}      & \text{Number of elements} \\
\text{number}    & \text{Frattini subgroup} & \text{of order 2}         \\ 
\text{of $G$}    & \text{of $V(KG)$}        & \text{in $V(KG)$}         \\ \hline               
3                & 256                      & 5119                      \\ \hline
4                & 256                      & 3^2 \cdot 5 \cdot 7 \cdot 13 \\ \hline
6                & 512                      & 5 \cdot 307               \\ \hline
11               & 128                      & 6143                      \\ \hline
12               & 128                      & 3^2 \cdot 5 \cdot 7 \cdot 13 \\ \hline
13               & 128                      & 3583                      \\ \hline
\end{array}
$$

It is easy to see from the table that these two invariants split all groups.

Note that generators of $V(KG)$ for 2-groups of orders up to 16
were given by R.~Sandling in \cite{Sandling}, but, as it was
confirmed by R.~Sandling in private communication, he did not
consider the isomorphism question.

\section{Groups of order 32}
\label{order32}

Among 51 groups of order 32, there are 
7 abelian groups and 3 groups of maximal class, so it remains to consider
41 groups of this order.

First we compute the order of the center of $V(KG)$
and the order of the Frattini subgroup of $V(KG)$ to split these groups
into twelve families, and then we will show that unit groups within each 
family are pairwise non-isomorphic. In the next table we list catalogue
numbers of groups and values of invariants for each family:

$$
\begin{array}{|c|c|c|c|c|}\hline
\text{Number} & \text{Catalogue} & \text{Order of the} & \text{Order of the} \\
\text{of the} & \text{numbers of} & \text{center} & \text{Frattini subgroup} \\ 
\text{family} & \text{groups of order 32} & \text{of $V(KG)$} & \text{of $V(KG)$} \\ 
\hline  
1  & 43,44                           & 2^{10} & 2^{23}   \\ \hline 
2  & 6,7,8                           & 2^{10} & 2^{24}   \\ \hline 
3  & 28,29,39,40,41,42               & 2^{13} & 2^{20}   \\ \hline 
4  & 9,10,13,14,27,30,31,32,33,34,35 & 2^{13} & 2^{21}   \\ \hline 
5  & 11,15                           & 2^{13} & 2^{22}   \\ \hline
6  & 49,50                           & 2^{16} & 2^{16}   \\ \hline 
7  & 46,47,48                        & 2^{19} & 2^{13}   \\ \hline 
8  & 22,23                           & 2^{19} & 2^{14}   \\ \hline
9  & 2,24,25,26,37,38                & 2^{19} & 2^{15}   \\ \hline 
10 & 5,12                            & 2^{19} & 2^{17}   \\ \hline 
11 & 4                               & 2^{19} & 2^{18}   \\ \hline 
12 & 17                              & 2^{19} & 2^{19}   \\ \hline 
\end{array}
$$

We see that families 11 and 12 contain only one group, so for these
groups the test is already finished. Now we will consider each family
1--10 separately. 

\noindent $\bullet$ Family 1. Since the number of involutions
in $V(KG)$ for both groups are the same, we need to employ other 
invariants. Using the \textsf{GAP4} Package \textsf{AutPGrp} by
B.~Eick and E.~O'Brien \cite{AutPGrp}, we computed the automorphism
groups $Aut(V(KG))$ for both cases. It appears that their orders,
given in the next table, are different:
$$
\begin{array}{|c|c|}\hline
\text{Catalogue}     & \text{Order of the automorphism}    \\
\text{number of $G$} & \text{group of $V(KG)$}             \\  \hline 
43 & \; 2^{102} \; \\
44 & \; 2^{101} \; \\ \hline
\end{array}
$$

\noindent $\bullet$ Family 2. The number of involutions in $V(KG_n)$ is
equal to $2^{18} \cdot 19$, $2^{19} \cdot 7$ and $2^{20} \cdot 3$ 
for $n=6,7$ and $8$, respectively. 

\noindent $\bullet$ Family 3. First we note that the exponent of the
center of $V(KG)$ is equal to 2 for $n \in \{28,29\}$ and to 4 for
$n \in \{39,40,41,42\}$. In case of exponent 2, automorphism groups
of $V(KG)$ have different orders, given in the following table:
$$
\begin{array}{|c|c|}\hline
\text{Catalogue}     & \text{Order of the automorphism}    \\
\text{number of $G$} & \text{group of $V(KG)$}             \\  \hline 
28 & \; 2^{175} \; \\
29 & \; 2^{173} \; \\ \hline
\end{array}
$$
In the case of exponent 4, the number of involutions in the 
center of $V(KG)$ is equal to 4095 for $n \in \{39,40,41\}$ and to
2047 for $n=42$, so the latter group is determined. Now it remains to 
check that the number of involutions in $V(KG_n)$ is equal to 
$2^{23}$, $2^{18} \cdot 31$ and $2^{19} \cdot 3 \cdot 5$ 
for $n=39,40$ and $41$, respectively. 

\noindent $\bullet$ Family 4. First we give a table containing
the number of involutions in $V(KG)$:
$$
\begin{array}{|c|c|}\hline
\text{Catalogue}     & \text{Number of elements}    \\
\text{number of $G$} & \text{of order 2 in $V(KG)$} \\  \hline  
9                    & 2^{18} \cdot 29              \\ \hline
10                   & 2^{20} \cdot 7               \\ \hline
13                   & 2^{19} \cdot 13              \\ \hline
14                   & 2^{19} \cdot 13              \\ \hline
27                   & 2^{18} \cdot 97              \\ \hline
30                   & 2^{18} \cdot 3^{4}           \\ \hline
31                   & 2^{19} \cdot 3 \cdot 11      \\ \hline
32                   & 2^{24}                       \\ \hline
33                   & 2^{18} \cdot 73              \\ \hline
34                   & 2^{20} \cdot 17              \\ \hline
35                   & 2^{24}                       \\ \hline
\end{array}
$$
We see that this parameter splits almost all groups 
except two pairs $n=13,14$ and $n=32,35$. To split these pairs, 
first we compute the order of the automorphism group of $V(KG)$:
$$
\begin{array}{|c|c|}\hline
\text{Catalogue}     & \text{Order of the automorphism}    \\
\text{number of $G$} & \text{group of $V(KG)$}             \\  \hline 
13 & \; 2^{149} \; \\
14 & \; 2^{149} \; \\
32 & \; 2^{158} \; \\
35 & \; 2^{159} \; \\ \hline
\end{array}
$$
Now to split the remaining pair for $n=13$ and $14$, we need one more step.
We were trying to compare a number of invariants but they were the same.
Finally, using the development version of the \textsf{AutPGrp} package
\cite{AutPGrp} we checked that the minimal number of generators
of $Aut(V(KG))$ is equal to 15 for $n=13$ and to 16 for $n=14$. 
Note that these groups, given by relations
$$
\begin{matrix}
G_{13} & = & \langle a,b \mid a^8=1, b^4=1 ,b^{-1}ab=a^3 \rangle ; \\
G_{14} & = & \langle a,b \mid a^8=1, b^4=1, b^{-1}ab=a^{-1} \rangle,
\end{matrix}
$$
also appeared in the preprint by Newman, Michler and O'Brien \cite{MNO},
where this pair also was the most complicated case in confirming {\bf(MIP)}
for groups of order 32.

\noindent $\bullet$ Family 5. The number of involutions in $V(KG_n)$ is
equal to $2^{20} \cdot 3$ for $n=11$ and $2^{18} \cdot 3^{2}$ for $n=15$.

\noindent $\bullet$ Family 6.
The number of involutions in $V(KG_n)$ is
equal to $2^{21} \cdot 7$ for $n=49$ and $2^{22} \cdot 3$ for $n=50$.

\noindent $\bullet$ Family 7. First we note that the exponent of the center
of $V(KG)$ is equal to 4 for $n=48$, while it is equal to 2 for $n=46,47$.
Then, the number of involutions in $V(KG_n)$ is equal to $2^{25} \cdot 3$ for 
$n=46$ and $2^{23} \cdot 11$ for $n=47$.

\noindent $\bullet$ Family 8. The number of involutions in $V(KG_n)$ is
equal to $2^{22} \cdot 23$ for $n=22$ and $2^{23} \cdot 11$ for $n=23$.

\noindent $\bullet$ Family 9. First we note that the exponent of the center
of $V(KG)$ is equal to 2 for $n=2$ and to 8 for $n=38$, while it is equal
to 4 for $n \in \{24,25,26,37\}$. Furthermore, the $p$-class of $V(KG)$ is
equal to 2 for $n \in \{2,24,25,26\}$ and to 3 for $n \in \{37,38\}$. Combining
these two invariants, we split groups with $n \in \{2,37,38\}$. For the three
remaining groups the number of involutions in $V(KG_n)$ is
equal to $2^{21} \cdot 3 \cdot 5$, $2^{22} \cdot 7$ and $2^{24}$ 
for $n=24,25$ and $26$, respectively. 

\noindent $\bullet$ Family 10. The number of involutions in $V(KG_n)$ is
equal to $2^{23} \cdot 11$ for $n=5$ and $2^{24}$ for $n=12$.

Thus, we confirmed {\bf UMIP} for groups of order 32, because each group 
$G$ of order 32 can be determined from the unique set of invariants of the
normalized unit group $V(KG)$ of its modular group algebra over the field
of two elements.

\section{Computational issues}
\label{comp}

Computations of $V(KG)$ and its invariants were performed using the
computational algebra system \textsf{GAP} \cite{GAP} enhanced with
the \textsf{LAGUNA} package \cite{LAGUNA}. Note that without the \textsf{LAGUNA}
package it would not be possible to compute the unit group of the modular group
algebra for groups of order 16 and 32.

First we tried to employ some "cheap" invariants such as
the order and the exponent of the center of $V(KG)$, 
and the Frattini subgroup of $V(KG)$. Computations of 
such kind could be performed very fast, especially on 
modern computers. 

For those groups, which can not be split using "cheap"
invariants, we decided to compare the number of involutions
in $V(KG)$. This parameter proved to be useful to deal with
{\bf (UMIP)} for 2-groups of maximal class \cite{BB1}, where
it was computed by theoretical means. We used the results of 
\cite{BB1} to check that our program in \textsf{GAP} returns
correct output for 2-groups of maximal class. We hope that 
results of our computations may motivate further attempts for
evaluation the number of involutions in unit groups and
provide a number of examples for such research.

For groups of order 16 it takes a very short time to count 
involutions in $V(KG)$, but this is already not the case for 
groups of order 32. This is why computations of the number 
of involutions for groups of order 32 were performed on the 
Computational Cluster of the Kiev National Taras Shevchenko 
University, created with the support of Intel 
corporation (\verb+http://www.cluster.kiev.ua/+), where it 
takes about 24 hours to enumerate elements of one unit group. 
Also, because of some internal \textsf{GAP} limitations, we
were forced first to enumerate cosets of $V(KG)$ by one of its
proper subgroup, and then to enumerate elements of these cosets.

In that cases when the number of involutions was not useful,
we tried to compare the number of elements of each order,
but this did not give any new information. It was the 
\textsf{AutPGrp} package \cite{AutPGrp} that brings new
information useful for splitting the most difficult pairs.
We especially acknowledge the help of B.~Eick in splitting 
the pair $(G_{13},G_{14})$ using the development version of 
\textsf{AutPGrp}, and for her helpful advises in \textsf{AutPGrp}
usage. Note that some other pairs of groups, which were split by
the number of involutions, could be also splitted in shorter time
by the order of $Aut(V(KG))$. Since this is easier for a \textsf{GAP}
user, we did not give these parameters here, being more interested in
the number of involutions by the above mentioned motivation.

\medskip
\noindent
\begin{minipage}[t]{7.5cm}
\small{Alexander Konovalov\\
       School of Computer Science\\
       University of St Andrews\\
       Jack Cole Building, North Haugh, \\
       St Andrews, Fife, KY16 9SX, Scotland\\
       e-mail: konovalov@member.ams.org}
\end{minipage}
\begin{minipage}[t]{7.5cm}
\small{Anastasiya Krivokhata\\
       Department of Mathematics\\
       Zaporozhye National University\\
       ul.Zhukovskogo, 66, Zaporozhye\\
       69063, Ukraine\\
       e-mail: k-algebra@zsu.zp.ua}
\end{minipage}
\end{document}